\documentclass[12pt,A4paper]{article}
\usepackage[T2A]{fontenc}
\usepackage[cp1251]{inputenc}
\usepackage[english]{babel}
\usepackage{amssymb}
\usepackage{cite}
\usepackage{latexsym}
\usepackage{amsmath}

\begin{document}

\medskip
\centerline{\bfseries ON PARTIALLY CONJUGATE-PERMUTABLE SUBGROUPS}
\centerline{\bfseries OF FINITE GROUPS}
\bigskip

\centerline{ V.I. Murashka}

\centerline{E-mail: mvimath@yandex.ru}

\centerline{and}

\centerline{A.F.Vasil'ev}

\centerline{E-mail:formation56@mail.ru}

\medskip
       \centerline{\textbf{Abstract.}}

       Let $ R $ be a subset of a group $ G $. We call a subgroup $ H $ of $ G $
 the $ R $-conjugate-permutable subgroup of $ G $, if $ HH^{x} = H^{x} H $ for all
 $ x \in R $. This concept is a generalization of conjugate-permutable subgroups
 introduced by T. Foguel.
Our work focuses on the influence of $R$-conjugate-permutable subgroups
 on the structure of finite groups in case when $R$ is the Fitting subgroup or its
 generalizations $F^{*}(G)$ (introduced by H. Bender in 1970) and $\tilde{F} (G)$
 (introduced by P. Shmid 1972).
We obtain a new criteria for nilpotency and supersoulbility of finite
 groups which generalize some well known results.

         \medskip

  \section {\bf Introduction} \label {g1}

All groups considered here are finite. Recall \cite{12} that the subgroups $ H $ and $ K $
of $ G $ are said to permute if $ HK = KH $, which is equivalent to that
set $ HK $ is a subgroup of $ G $.

The classic area of group theory is the study of subgroups of $ G $ which permute with every
 subgroup of a dedicated system of subgroups of $ G $. This trend goes back to O. Ore \cite{O} who
introduced the concept of quasinormal (permutable) subgroup in 1939. Recall that
subgroup $ H $ of a group $ G $ is called quasinormal if it permutes with each subgroup of $ G $.
 Every normal subgroup is quasinormal. It is known that every quasinormal subgroup
is subnormal. There are examples showing that the converse
need not hold.

Another important type of subgroups' permutability was proposed by
O. Kegel  \cite{K} in 1962. A subgroup $ H $ of a group $ G $ is called
$ S $-permutable ($ S $-quasinormal, $ \pi $-quasinormal)
subgroup of $ G $, if $ H $ permutes with every Sylow subgroup of
$ G $. Note that every $ S $-permutable subgroup is
subnormal. The converse need not hold. Currently, the concept of quasinormal and
$ S $-permutable subgroups and their generalizations have been studied intensively by many authors
(see reviews \cite{4}, \cite{5} and the monograph \cite{6}).

In 1997 T. Foguel \cite{7} noted in the proof that a quasinormanal subgroup is subnormal, one only needs
 to show that it is permutable with all of its conjugates.
 This led him to the following  concept of subgroups' permutability.

\medskip

 {\bf Definition 1 \cite{7}.} {\sl A subgroup $ H $ of a group $ G $ is called the conjugate-permutable
 subgroup of $ G $, if $ HH ^ {x} = H ^ {x} H $ for all $ x \in G $. Denoted by
$ H <_ {\texttt {C-P}} G $}.

\medskip

Clearly, every quasinormal subgroup is conjugate-permutable. In \cite{7} there is an example
 showing that the converse is not true. On the other hand, every 2-subnormal subgroup
 (i.e. subgroup is a normal subgroup of some normal subgroup of the group) is a conjugate-permutable.

The concept of conjugate-permutable subgroups proved useful. A number of
authors (for example see \cite{7},\cite{8},\cite{9},\cite{10},\cite{11}) studied the influence of conjugate-permutable subgroups on the structure of the group.
Particulary, some criteria of nilpotency of finite groups were obtained in
 terms of the conjugate-permutable subgroups  (see \cite{7},\cite{9}).

Analyzing the proofs of some results of the works \cite{7},\cite{9}, we have seen that
a smaller number of conjugate subgroups $ H^{x} $ which $ H $ must permute with
 can be considered  for the conjugate-permutable subgroup $ H $.
This observation led us to the following definition.

\medskip

  {\bf Definition 2.} {\sl Let $ R $ be a subset of a group $ G $.  We shall call a subgroup
 $ H $ of $ G $ the  $R $-conjugate-permutable subgroup of $ G $, if $ HH^{x} = H^{x} H $
 for all $ x \in R $. Denoted by $ H <_{R-C-P} G $.}

\medskip

Note that if $ R = G $, then $ H $ is a conjugate-permutable, if $ R = 1 $, then
any subgroup of $ G $ is $ R $-conjugate-permutable subgroup.

In this paper we study the influence of various systems of $ R $-conjugate-permutable subgroups on
  the structure of  finite groups $G$, when $ R \in \{F(G), F^{*} (G), \tilde {F}(G) \} $. Recall that

 The Fitting subgroup $ F (G) $ is maximal normal nilpotent subgroup of $ G $.

 $ \tilde {F}(G) $ is a generalization of the Fitting subgroup introduced by P. Schmid (\cite{13}, p. 79,
 or \cite{17}]. It is defined by $  \Phi (G) \subseteq \tilde {F}(G) $ and
 $  \tilde {F} (G) / \Phi (G) = Soc (G / \Phi (G)) $.

A subgroup $ F^{*} (G) $ is another generalization of the Fitting subgroup.
 It was introduced by H. Bender \cite{18} and defined by $ F^{*} (G) / F (G) = Soc (C_{G} (F (G) F (G) / F (G)) $.

As we will show that $F(G) \subseteq F^{*} (G)\subseteq \tilde {F}(G) $ for every group $G$.
 If $ G $ is a solvable group, then $ \tilde {F}(G)=F^{*} (G)=F(G) $.
\medskip

{\bf Example 1.} Let $ G \simeq S_{4} $ be the symmetric group of degree 4. Let $ H $ be
 Sylow 2-subgroup of $ G $. Then $ H $ is a maximal subgroup of $ G $ which
 is not normal in $ G $. Note that $ \tilde {F} (G)=F^{*}(G) = F (G) \subseteq H $.
 Hence $H<_{F(G)-C-P}G$. By the theorem of Ore (\cite{12} p. 57)
 $ HH^{x} \neq H^{x}H $ for all $ x \in G \backslash H $. Therefore $ H $ is
 not conjugate-permutable subgroup of $ G $.  
\medskip

{\bf Theorem A.} {\sl A group $ G $ is nilpotent if and only if every maximal subgroup of $G$ is
 $ \tilde {F} (G) $-conjugate-permutable.}

The following example shows that we can not use $F^{*}(G)$ in place of
$\tilde{F}(G)$ in theorem A.

{\bf Example 2.} Let $ G \simeq A_{5}$ be the alternating group of degree
5, and $ K = F_{3}$ be a field composed by three elements. We denoted by
$A = A_{K}(G)$ the Frattini $KG$-module \cite{20}. In view of
\cite{20}, $A$ is  an irreducible $KG$-module of the dimension 4. By known Gaschutz theorem, there
exists a Frattini extension  $A\rightarrowtail R\twoheadrightarrow G$
such that $A\stackrel {G}{\simeq} \Phi(R)$ and $R/\Phi(R)\simeq G$.
From the properties of module $A$ follows that $\tilde {F}(G) = R$ and  $F^{*} (G) = \Phi(R)$.
Note that every  maximal subgroup of $G$ is  $ F^{*}(G)$-conjugate-permutable, but the group $R$
 is not nilpotent.

   {\bf Corollary A.1 \cite{7}}. {\sl If every maximal subgroup a group
$G$ is conjugate-permutable then $G$ is nilpotent.}

{\bf Corollary A.2}. {\sl If $G$ is a non-nilpotent group then there is an abnormal  maximal subgroup $M$ of $G$
 such that $\tilde{F}(G)\nsubseteq M$.}

\medskip

 The reference system of subgroups of a
group is the set of its Sylow subgroups, the knowledge of the structure
and embedding properties allows in many cases to reveal the structure of the group.
For example, recall the following well-known result: the group is nilpotent
if and only if each of its Sylow subgroup is subnormal.

\medskip

{\bf Theorem B.} {\sl The following statements for a group $G$ are equivalent:

1) $ G $ is nilpotent;

2) every abnormal subgroup of $ G $ is $ F^{*}(G) $-conjugate-permutable
 subgroup of $ G $;

3) normalizers of all Sylow subgroups of $ G $ are
 $ F^{*} (G) $-conjugate-permutable subgroups of $ G $;

4)  Sylow subgroups of $ G $ are
$ F^{*} (G) $-conjugate-permutable subgroups of group $ G $.}

{\bf Corollary B.1.} {\sl A group $ G $ is nilpotent if and only if
 the normalizers of all Sylow subgroups of $ G $ contains
$ F^{*}(G) $}.

{\bf Corollary B.2.} {\sl A group $ G $ is nilpotent if and only if it satisfies one of the
 following conditions:

  1) Sylow subgroups of maximal subgroups of $ G $ are $ F^{*}(G) $-conjugate-permutable;

  2) normalizers of Sylow subgroups in maximal subgroups of $G$ are $ F^{*}(G) $-conjugate-permutable.}

{\bf Theorem C.} {\sl If all cyclic primary subgroups of
 $ G $ are $ F^{*}(G) $-conjugate-permutable then $ G $ is nilpotent.}

As follows from example 1.2 \cite{7}  the converse of theorem C is false.

{\bf Corollary C.1.} {\sl A group $ G $ is nilpotent if  cyclic subgroups of Sylow subgroups of
 maximal subgroups of $G$ are $ F^{*}(G) $-conjugate-permutable; }

\medskip

Note that theorems 2.3, 2.4 \cite{9} and theorem 2.3 \cite{7} follows from theorems  A, B and C.

    A group $ G $ is called dinilpotent [6, p. 100] if $ G = AB $, where $ A $ and $ B $ are nilpotent
subgroups of $ G $. By the theorem of Wielandt-Kegel every
dinilpotent group is solvable. Examples of dinilpotent groups are
biprimary, supersoluble groups and other. Dinilpotent group
studied by many authors in different directions (see \cite{6}, Chapter 3). In this
direction we obtain the following result.

\medskip

{\bf Theorem D.} {\sl Let $ G = AB $ be a dinilpotent group. The group $ G $ is nilpotent
 if and only if $ G = AB $, where $ A $ and $ B $ are $ F(G) $-conjugate-permutable
 subgroups of $ G $.}

 {\bf Corollary D.1} {\sl If the group $ G = AB $ is a dinilpotent group and
$ F (G) \subseteq A \cap B $ then $ G $ is nilpotent.}

 It is well known that the product of two normal supersoluble subgroups of
a group is not necessarily supersoluble. This fact has been the starting point for a series of results about
factorized groups in which the factors satisfy certain permutability
conditions. The first step in this direction was taken by Baer \cite{21}.
He proved that if $G=AB$ is the product of the supersoluble normal
subgroups $A$ and $B$ and $G'$ is nilpotent then $G$ is supersoluble. We
obtain the following result.

\medskip

{\bf Theorem E.} {\sl A group $ G $ is supersoluble if and only if
$ G = AB $, where $ A $ and $ B $ are supersoluble $F(G)$-conjugate-permutable
 subgroups of $ G $ and $G'$ is nilpotent.}

\medskip

It is well known ( \cite{22}, p. 127) that if a group $G$ contains two normal supersoluble subgroups with coprime 
indexes in $G$ then $G$ is supersoluble. For $F(G)$-conjugate-permutable subgroups this result is not true. 

\textbf{Example 3.} Let $G$ be the symmetric group of degree 3. By theorem 10.3B \cite{12} there is a faithful
 irreducible $F_{7}G$-module $V$ over the field $F_{7}$ of 7 elements and the dimension of $V$ is 2.
 Let $R$ be the semidirect product of $ V$ and $G$. Let $A=VG_{3}$ and $B=VG_{2}$ where $G_{p}$ is a
 Sylow $p$-subgroup of $G$, $p\in{2,3}$. Since $7\equiv 1 ($mod$p)$ for $p\in{2,3}$, it is easy to
 check that subgroups $A$ and $B$ are supersoluble. Since $V$ is faithful irreducible module, $F(R)=V$.
 Therefore $A$ and $B$ are the $F(R)$-conjugate-permutable subgroups of $G$. Note that $R=AB$ but $R$ 
 is not supersoluble.

\textbf{Theorem F.} {\sl A group $G$ is supersoluble if and only if $G$ is metanilpotent and contains 
 two supersoluble $F(G)$-conjugate-permutable subgroups with coprime 
indexes.}

\textbf{Theorem G.} {\sl A group $G$ is supersoluble if and only if $G=AB$ is the product of
 supersoluble $F(G)$-conjugate-permutable subgroups $A$ and $B$ and 
$G' = A'B'$.}

\medskip

\section {\bf Preliminaries} \label {g2}

\medskip

We use standard notation and terminology, which if necessary can be found
 in \cite{6}, \cite{12}, \cite{13} and \cite{21}. We recall some definitions, notions and results.

 Let $G$ be a group then $ Syl_{p} (G) $ denote the set of all Sylow
 $ p $-subgroups of $ G $; $ H \leq G $ means $ H $ is a subgroup of $ G $; $ | G | $ is the order of
  $ G $; $ H \triangleleft G $ means $ H $ is a normal subgroup of $ G $;
 $ H \triangleleft \triangleleft G $ means $ H $ is a subnormal subgroup of $ G $;
 $ H_ {G} $ is the core of the subgroup $ H $ of $ G $, i.e. maximal normal subgroup of $G$ contained in
 $ H $; $ C_{G} (H) $ is the centralizer of the subgroup $ H $ in $ G $; $ N_{G} (H) $
 is the normalizer of subgroup $ H $ in $ G $; $ Z (G) $ is the center of $ G $;
 $ Z_{\infty} (G) $ is the hypercenter of
 $ G $; $ Soc (G) $ is socle of $ G $, i.e. the product of all minimal normal subgroups of
 $ G $. We will use the symbol 1 to denote the identity subgroup  of a 
group.

We recall the following well-known definitions and results (see  \cite{12}):

{\bf Definition 2.1.} {\sl The group is called nilpotent if all its
Sylow subgroups are normal.}

{\bf Theorem 2.2.} {\sl For a group $ G $ the following statements
are equivalent:

A) $ G $ is nilpotent;

2) $ G $ is a direct product of its Sylow subgroups;

3) every proper subgroup of $G$ is distinct from its normalizer;

4) all maximal subgroups of $G$ are normal;

5) all subgroups of $ G $ are subnormal.}

{\bf Definition 2.3.} {\sl A  subgroup $ M $ of  non-unit group
$ G $ is the  maximal subgroup of $G$ if $ M $ is not contained in any
other subgroup distinct from  $ G $.}

{\bf Definition 2.4.} {\sl  The Frattini subgroup of a $G \neq 1$ is  the
intersection of all  maximal subgroups of $G$. Denoted by $ \Phi (G) $}.

 Set that $ \Phi (1) = 1 $.

{\bf Theorem 2.5.} {\sl $ \Phi (G) $ is a normal
nilpotent subgroup of $G$}.

{\bf Theorem 2.6.} {\sl Let $ D \triangleleft K \triangleleft G $, $ D \leq \Phi (G) $
and $ D \triangleleft G $. If the quotient group $ K / D $ is
nilpotent then  $ K $ is nilpotent.}

\textbf {Lemma 2.7 (Frattini)}. {\sl If $ K $ is a normal subgroup of a group $ G $ and $ P $
 is a Sylow subgroup of $ K $ then $ G = N_{G} (G) K $.}

\textbf {Definition 2.8.} {\sl Let $G$ be a group. A subgroup $ H $ is called pronormal in $ G
$ if
 subgroups $ H $ and $ H^{x} $ are conjugate in $ \langle H, H^{x}
\rangle $ for all $x \in G $.}

\textbf {Definition 2.9.} {\sl Let $G$ be a group. A subgroup $ H $ is called abnormal of $ G $ if
  $ x \in \langle H, H^{x} \rangle $ for all $ x \in G $.}

\textbf {Lemma 2.10.} {\sl If a subgroup $ H $ is pronormal in $ G $ then
$ N_{G} (H) $ is abnormal subgroup of $ G $.}

\textbf {Lemma 2.11.} {\sl Let $ H $ be an abnormal subgroup of a group $ G $.
 From $ H \leq U \leq G $ and $ H \leq U \bigcap U^{x} $ follows that $ x
\in U $.}

\textbf {Theorem 2.12.} {\sl Let $ \Phi (G) = E $. Then
$ Z_{\infty} (G) = Z (G) $.}

\textbf {Theorem 2.13 \cite{15}.} {\sl The hypercenter  is the intersection of all the normalizers
 of all Sylow subgroups.}

\textbf {Theorem 2.14}. {\sl Let $ G $ be a group. The Fitting subgroup  $ F (G / \Phi (G)) = F (G) / \Phi (G) $
and is equal to the direct product of abelian minimal normal subgroups of $ G / \Phi (G) $.}

The idea of lemma's 2.15 proof was proposed by L. Shemetkov.

\textbf{Lemma 2.15.} {\sl
$F^{*}(G)\subseteq \tilde {F}(G)$ for any group $G$.}

\emph{Proof}. Let a group $G$ be the minimal order counterexample for
lemma 2.15.
If $\Phi(G)\neq 1$ then for $G/\Phi(G)$ the statement is true. From
 $F^{*}(G)/\Phi(G)\subseteq F^{*}(G/\Phi(G))$ and  $\tilde{F}(G/\Phi(G))=\tilde{F}(G) /\Phi(G)$
we have that $F^{*}(G)\subseteq \tilde {F}(G)$. It is a contradiction with the choice of $G$.

Let $\Phi(G) = 1$. Now $\tilde{F}(G)= Soc(G)$. By 13.14.X  \cite{19}  $F^{*}(G) = E(G)F(G)$.
 Note $\Phi(E(G)) = 1$. Since 13.7.X  \cite{19}  $E(G)/Z(E(G))$ is the direct product of simple
 nonabelian groups, $Z(E(G)) = F(E(G))$. From it and theorem 10.6.A.c \cite{12} we conclude that
 $E(G) = HZ(E(G))$ where $H$ is the complement to $Z(E(G))$ in $E(G)$. Now $H$ is the direct product
 of simple nonabelian groups. Since $H char E(G)\triangleleft G$, we have $H\triangleleft G$. From
 lemma 14.14.A \cite{12} follows $H\subseteq Soc(G)$. Since  $Z(E(G))\subseteq F(G)\subseteq \tilde{F}(G)$
   and $H\subseteq Soc(G)$, it follows that
  $E(G)\subseteq \tilde{F}(G)$. Now $F^{*}(G) = E(G)F(G)\subseteq \tilde{F}(G)$.
 It is a contradiction with the choice of $G$. Lemma is proved.

The subgroups $ \tilde{F} (G)$, $F^{*}(G)$ and $F(G)$ have the following useful properties
 \cite{12}, \cite{13}, \cite{19}.

\textbf{Lemma 2.16.} {\sl Let $G$ be a group. Then

1) $\tilde{F}(G/\Phi(G))=\tilde{F}(G)/\Phi(G)$;

2) $C_{G}(\tilde{F} (G))\subseteq F(G)$.

3) if $G$ is soluble  $C_{G}(F(G))\subseteq F(G)$;

4) $C_{G}(F^{*}(G))\subseteq  F(G)$;

5) $F(G)\subseteq F^{*}(G)\subseteq \tilde{F} (G)$}.

{\bf Lemma 2.17.} {\sl Let $G = AB$ be a product of the normal nilpotent subgroup $A$ and
 subnormal supersoluble subgroup $B$. Then $G$ is supersoluble}.

\emph {Proof.}    Follows from theorem 15.10 \cite{13}.

{\bf Definition 2.18}. {\sl A group $G$ is called quasinilpotent if $G=C_{G}(H/K)H$ for any chief
 factor $H/K$ of $G$.}

{\bf Theorem 2.19}. {\sl A group $G$ is quasinilpotent if and only if 
$G/Z_{\infty}(G)$ is quasinilpotent.}

Follows from theorem 13.6 \cite{19}.

 {\bf Definition 2.20}. {\sl A group is said to be supersoluble whenever 
its chief factors are all cyclic.}

 {\bf Definition 2.21}. {\sl The commutator subgroup of a supersoluble 
group is nilpotent}.

\textbf{Theorem 2.22.} {\sl If $G$ is the product of two subnormal supersoluble subgroups with coprime
 indexes in $G$ then $G$ is supersoluble.}

Theorem 2.20 follows from theorem 3.4 (p. 127 \cite{22}) and the induction by the order of $G$.

     \medskip

          \section {\bf Properties of $ R $-conjugate-permutable Subgroups} \label {g3}

\medskip

{\bf Lemma 3.1} \cite{7}.  {\sl Let $H<_{C-P}G$ then $H \triangleleft
\triangleleft G$}.

{\bf Lemma 3.2.} {\sl Let $ H <_ {R-C-P} G $, $ R \leq G $ and $ RH = HR $. Then
 $ H \triangleleft \triangleleft HR $.}

\emph {Proof.} Let $ k \in HR $. Then $ k = hr $, where $ h \in H $ and $ r \in R $.
 Since $ H <_{R-C-P} G $, $ HH ^ {k} = HH ^ {hr} = HH ^ {r} = H ^ {r} H = H ^ {hr} H = H ^ {k} H $.
 Therefore $ H <_{C-P} HR $. Now the result follows from lemma 3.1. Lemma is proved.

{\bf Lemma 3.3.} {\sl Let $ H <_{C-P} G $ and $ H $ is pronormal in $ G $. Then $ H \triangleleft G $}.

\emph {Proof}. Since $ H <_{C-P} $, then by lemma 3.1 $ H \triangleleft \triangleleft G $.
Consider  $ N_{G} (H) $. By Lemma 2.10 $ N_{G} (H) $ is abnormal subgroup of
$ G $.
Since $ H \triangleleft \triangleleft G $, $ N_{G} (H) \triangleleft \triangleleft G $.
Assume that $ N_{G} (H) \neq G $. Then there exists $ x $ not in $ N_{G} (H) $ such that $ N_{G} (H)^{x} = N_{G} (H) $. Then $ x $ does not belong
to $ \langle N_{G} (H)^{x} ; N_{G} (H) \rangle$. This contradicts the fact that
$ N_{G} (H) $ is abnormal subgroup of $ G $. It remains to assume that $ N_{G} (H) = G $ and
 $ H \triangleleft G $.
 The lemma is proved.

{\bf Proposition 3.4.} {\sl Let $ H <_{R-C-P} G $, $ R \leq G $ and $ RH = HR $, then
 $ H \triangleleft \triangleleft HR $.
 In particular, if $ H $ pronormal subgroup of $ HR $, then $ H \triangleleft HR $.}

\emph {Proof}. This follows from Lemmas 3.2 and 3.3.

\textbf {Corollary 3.4.1}. {\sl Let $ P <_{R-C-P} G $, $ P \in Syl_ {p} G $, $ R \leq G $ and
 $ PR = RP $.
 Then $ P \triangleleft PR $.}

\textbf {Corollary 3.4.2}. {\sl Let $ M <_{R-C-P} G $, $ R \leq G $, $ MR = RM $ and $ M $
 is maximal in $ MR $.
 Then $ M \triangleleft PR $.}

\textbf {Lemma 3.5.} {\sl Let $ H <_{C-P} G $ and $ K \triangleleft G $. Then $ HK <_{C-P} G $}.

 \medskip

 \section {\bf Results} \label {g4}

 \medskip

 {\bf Proof of Theorem A.} Let $ G $ be a nilpotent group. Then $ \tilde {F} (G) = G $. By 4)
Theorem 2.2 every maximal subgroup of $ G $ is normal in
$ G $, hence, $ \tilde {F} (G) $-conjugate-permutable.

  Conversely. Assume the result is false and $ G $ be a counterexample of minimal order.
Then $ G $ is the non-nilpotent group and  all
 maximal subgroups of $ G $  are $ \tilde {F} (G) $-conjugate-permutable.

 Suppose that $ \Phi (G) \neq 1 $. Consider the quotient $ G / \Phi (G) $.
We have $ \tilde {F} (G / \Phi (G)) = \tilde {F} (G) / \Phi (G) $. It is easily seen
 that all  maximal
 subgroups of $ G / \Phi (G) $ are $ \tilde {F} (G / \Phi (G)) $-conjugate-permutable.
 Since $ | G |> | G / \Phi (G) | $, we have $ G / \Phi (G) $ is nilpotent.
 From Theorem 2.6. follows that $G$ is nilpotent, a contradiction.

Assume that $ \Phi (G) = 1 $. Then   $ \tilde {F} (G) = Soc (G) $.

Assume now that $ \tilde {F} (G) $ is not nilpotent. By 2) theorem 2.2 there is a subgroup
 $ S \in Syl_ {p} (\tilde {F} (G)) $ such that $ S $ is not normal
 in $ \tilde {F} (G) $. Let $ P \in Syl_{p} (G) $ and $ P \cap \tilde {F} (G) = S $. Note
that $ S^{x} = P^{x} \cap \tilde {F} (G)^{x} = P \cap \tilde {F} (G) = S $ for every $ x \in N_{G} (P) $.
It means that $ N_{G} (P) \subseteq N_{G} (S) $.  Since $ N_{G} (S) \neq G $, we have $ N_{G} (P) \neq G $.
Let $ M $ be a maximal subgroup of $ G $  such that $ N_{G} (S)\subseteq M $.
By lemmas 2.10 and 2.11 $ M $ is the abnormal subgroup of $ G $. By lemma
2.7
 $ N_{G} (S) \tilde {F} (G) = M \tilde {F} (G) = G $. Since $ M $ is the
 $ \tilde {F} (G) $-conjugate-permutable subgroup, $ M $ is normal in $ G $ by Corollary 3.4.2,
 a contradiction.

Therefore we have that $ \tilde {F} (G) $ is nilpotent. Then
$ \tilde {F} (G) = F (G) = Soc (G) = N_{1} \times ... \times N_{t} $ where $ N_{i} $ runs over all
 minimal normal subgroup of $ G $. From $ \Phi (G) = 1 $ and Theorem 2.14,
 it follows that $ N_{i} $ is an abelian subgroup  for all $ i = 1, ..., t $. Then there is a maximal
 subgroup $ M_{i} $ such that $ N_{i} M_{i} = G $  for all $ i = 1, ..., t $.
 Note that $ M_{i} \tilde {F} (G) = G $.
 Since $ M_ {i} <_{\tilde {F} (G)-C-P}G $ , $ M_{i} $ is normal in $ G $
 for all $ i = 1,. .., t $ by Corollary 3.3.2.  Since $ N_{i} $ is abelian subgroup,
we have $ N_{i} \subseteq C_{G} (N_{i}) $ and $ N_{i} \cap M_ {i} = 1$. Then $ M_{i} \triangleleft G $
implies  $ M_{i} \subseteq C_{G} (N_{i}) $ for all $ i = 1, ..., t $. We
 show that $ G = M_ {i} N_ {i} \subseteq C_{G} (N_{i}) $ for every $ i = 1, ..., t $.
 Therefore $ N_{i} \subseteq Z (G) $ for all $ i = 1, ..., t $. Then
 $ \tilde {F} (G) \subseteq Z (G) $. Hence $ G \subseteq C_{G} (\tilde {F} (G)) \subseteq F (G) $.
 Thus $ G $ is nilpotent, a contradiction. Theorem A is proved.

\medskip

  {\bf Proof of Theorem B.} Prove that 1) implies 2). Since $G$ is nilpotent, $ F^{*} (G) = G $.
 Any subgroup of $ G $ is subnormal by 5) theorem 2.2. It means that the subgroup $G$ is the only
 one abnormal and subnormal subgroup in $G$. It is clear that $ G $ is the
 $ F^{*}(G) $-conjugate-permutable subgroup. Thus 1) implies 2).

 Normalizers of  all Sylow subgroups are abnormal subgroups by lemma 2.10. Therefore
2) implies 3).

 Assume the validity of 3). Let $ P $ be a Sylow subgroup of $ G $. Since $ N_{G} (P)$ is the
 $  F^{*}(G) $-conjugate-permutable subgroup and $ N_{G} (P) $ is pronormal
 in $ G $,  we have $ N_{G} (P) \triangleleft N_{G} (P) F^{*} (G) $  by lemmas 3.2
 and 3.3. By Lemma 2.10 $ N_{G} (P) $ is abnormal in $ G $. Therefore
$ N_{G} (P) $ is abnormal in  $ N_{G} (P) F^{*} (G) $. Then
 $ N_{G} (P) = N_{G} (P) F^{*} (G) $.
Which implies that $ P $ is normal in $ N_{G} (P) F^{*} (G) $. Hence $ P $ is
  the $ F^{*} (G) $-conjugate-permutable
subgroup of $ G $. Thus 3) implies 4).

Finally we show that 4) implies 1). Assume that 1) is not true and $G$ is a counterexample of least order.

Let $ P $ be  a Sylow subgroup of $ G $. Since $ P $ is the
 $ F^{*} (G) $-conjugate-permutable subgroup, we see that  $ P \triangleleft P F^{*} (G) $
 by Corollary 3.4.1.
It follows that $ F^{*} (G) \subseteq N_{G} (P) $.
 Now we have that $ F^{*} (G)$ lies in the intersection of all the normalizers
 of all Sylow subgroups of $ G $. Therefore $ F^{*} (G) \subseteq Z_{\infty} (G) $ by Theorem 2.13.
 Note that $F(G)=F^{*}(G)=Z_{\infty}(G)$.

Assume that $\Phi(G)\neq E$.  Let $H/\Phi(G)=F^{*}(G/\Phi(G))$. Show that $H/\Phi(G)=F^{*}(G)/\Phi(G)$.
 It is clear that
 $F^{*}(G)/\Phi(G)\subseteq H/\Phi(G)$. Suppose that
 $H/\Phi(G)/F^{*}(G)/\Phi(G)\simeq H/F^{*}(G)\neq 1$. Note that
 $H/\Phi(G)$ and $F^{*}(G)/\Phi(G)$ are quasinilpotent. It follows that $ H/F^{*}(G)$
 is quasinilpotent.
Now $H/Z_{\infty}(G)/Z_{\infty}(H)/Z_{\infty}(G)\simeq H/Z_{\infty}(H)$  is quasinilpotent.
 By theorem theorem 2.19 $H$ is the normal quasinilpotent subgroup of $G$. Hence $H\subseteq F^{*}(G)$.
 We have the contradiction with $H/F^{*}(G)\neq E$. Thus $F^{*}(G/\Phi(G))=F^{*}(G)/\Phi(G)$.

 Let $ S / \Phi (G) $ be a Sylow subgroup of $ G / \Phi (G) $. There is a Sylow subgroup $ P $ of $ G $
 such that $ P \Phi (G) / \Phi (G) = S / \Phi (G) $.
  From $ F^{*}(G / \Phi (G)) = F^{*} (G) / \Phi (G) $ it follows that $ S / \Phi (G) $ is the
 $ F^{*}(G / \Phi (G)) $-conjugate-permutable subgroup of $ G / \Phi (G) $. By minimality of $G$ we have
 that $ G / \Phi (G) $ is nilpotent. Hence $G$ is nilpotent by theorem 2.6, a contradiction.

Suppose now that $ \Phi (G) = 1 $.  By theorem 2.12 we have $ Z_{\infty} (G) = Z (G) $. Therefore
 $F^{*} (G) = Z (G) $.
 Now we have $ G = C_{G} (F^{*}  (G)) \subseteq F^{*} (G) $. Thus $ G $ is
 nilpotent. This is the final contradiction. Theorem B is proved.

 \medskip

{\bf Proof of Theorem C}.  Let $ P \in Syl_{p} (G) $ and $ x \in P $. Then the subgroup
 $ \langle x \rangle $ is the $ F^{*} (G) $-conjugate-permutable subgroup. By
lemmas 3.1 and 3.2
 $ \langle x \rangle \triangleleft \triangleleft \langle x \rangle F^{*} (G) $. Note
that $ \langle x \rangle \triangleleft \triangleleft P $. Since
$ \langle x \rangle \leq P \bigcap \langle x \rangle F^{*} (G) $ , by
theorem 1.1.7  (\cite{6} p.3) $ \langle x \rangle $ is the subnormal subgroup in the product
 $ P (\langle x \rangle F^{*} (G)) $.
 Since $ P $ is generated by its cyclic subnormal in $ P F^{*} (G) $ subgroups, by
  theorem (\cite{13} p. 70) we have that $ P \triangleleft \triangleleft PF^{*} (G) $. Then
$ P \triangleleft P F^{*} (G) $  by lemma 3.3. Thus $ F^{*} (G) \subseteq N_ {G} (P) $.
Now theorem C immediately follows from theorem B.

\medskip

  {\bf Proof of Theorem D}. Assume the result is not true and $ G $ is a minimal counterexample
 to  theorem D. Then $ G = AB$  is the non-nilpotent dinilpotent group
where  subgroups  $ A $ and $ B $ are $ F (G) $-conjugate-permutable.
 By well-known theorem of Wielandt-Kegel group $ G $ is solvable.

Assume that $ \Phi (G) \neq E $.  Since
 $ F (G / \Phi (G)) = F (G) / \Phi (G) $,  we have that $ A \Phi (G) / \Phi (G) $ and
 $ B \Phi (G) / \Phi (G) $ are
 the $ F (G / \Phi (G)) $-conjugate-permutable  nilpotent subgroups of $ G / \Phi (G) $. Now
 $ G / \Phi (G) = A \Phi (G) / \Phi (G) B \Phi ( G) / \Phi (G) $.
 By minimality of $G$ we have that $ G / \Phi (G) $ is nilpotent. Hence $G$ is nilpotent by
 theorem 2.6, a contradiction.

Assume that $ \Phi (G) = 1 $. Let $ R = AF (G) $. Then $ A \triangleleft \triangleleft R $ by
 lemma 3.2. Since  $ A $ and $ F (G) $ are nilpotent, we have that $ R $ is nilpotent. By analogy 
 $ T = BF (G) $
 is nilpotent. Since $ G $ is not nilpotent,
 there is an abnormal maximal  subgroup $ M $ of $ G $ such that $F(G)\nsubseteq M$ by Corollary A.2.

Note that $ G / M_{G} = RM_{G} / M_{G} TM_{G} / M_{G} $. By lemma 4 of \cite{14} without
 loss of generality we can assume that
 $ AM_{G} / M_{G} $ is $ p $-group, and
 $ BM_{G} / M_{G} $ is $ p '$-group. Now $ RM_{G} / M_{G} \bigcap TM_{G} / M_{G} = M_{G} / M_{G} $.
 It means $ RM_ {G} \bigcap TM_ {G} = M_ {G} $. Now
 it follows that $ F (G) \subseteq R \bigcap T \subseteq M_{G} $. But then
 $ F (G) \subseteq \bigcap M_{G} $,
where $ M$ runs  over all maximal subgroups $ M $ such that $ MF (G) \neq G $.
 By the  theorem of \cite{16} we have $ F (G) \subseteq \Phi (G) = 1 $. This is the final
 contradiction. Theorem D is proved.

 {\bf Proof of Theorem E}. Since $G'$ is nilpotent, $G$ is soluble. Since $F(G)$ is nilpotent, $K=AF(G)$
 is supersoluble by lemma 2.17. By analogy $R=BF(G)$ is supersoluble. Now $K/F(G) R/F(G)=G/F(G)$
 is abelian, since $G/G'$ is abelian and $G'$ is nilpotent. Therefore $K/F(G) $ and $ R/F(G)$ are normal
 in $G/F(G)$. Now $K$ and $R$ are normal in $G$. Now our theorem follows from Baer's theorem \cite{21}.

 {\bf Proof of Theorem F}. Since $G$ is metanilpotent, $G$ is soluble. Since $F(G)$ is nilpotent,
 $K=AF(G)$ is supersoluble by lemma 2.17. By analogy $R=BF(G)$ is supersoluble. Now $K/F(G) R/F(G)=G/F(G)$
 is nilpotent. By 3) theorem 2.2 $K/F(G)$ and $  R/F(G)$ are subnormal subgroups of $G/F(G)$. Therefore
 $K$ and $R$ are subnormal subgroups of $G$. By theorem 2.20 $G$ is supersoluble.

 {\bf Proof of Theorem G}. Since $F(G)$ is nilpotent, $K=AF(G)$ is supersoluble by lemma 2.17.
 By analogy $R=BF(G)$ is supersoluble. Since $A\leq K$ and $B\leq R$, 
  we see that $A'\leq K'$ and $B'\leq R'$.  By $G'= A'B'= R'K'$ is 
the  product of nilpotent subgroups $R'$ and $K'$ because the commutant of a supersoluble group is nilpotent.
 Note that $K'<_{F(G)-C-P}K$ and $R'<_{F(G)-C-P}R$. Since $G' \triangleleft G$, $F(G')\leq F(G)$.
 So $G'$ is the product of two nilpotent $F(G')$-conjugate-permutable subgroups. From theorem D it follows that
 $G'$ is nilpotent.
 This and theorem E immediately combine to yield.

\end{document}